\documentclass[a4paper,12pt]{amsart}
\usepackage{ifthen}
\usepackage{mathrsfs}
\nonstopmode
\numberwithin{equation}{section}
\newtheorem{thm}{Theorem}[section]
\newtheorem{cor}[thm]{Corollary}
\newtheorem{lem}[thm]{Lemma}
\newtheorem{prop}[thm]{Proposition}

\theoremstyle{definition}

\newenvironment{pf}[1][]{%
 \vskip 3mm
 \noindent
 \ifthenelse{\equal{#1}{}}%
  {{\slshape Proof. }}%
  {{\slshape #1.} }%
 }%
{\qed\bigskip}

\newcounter{alphabet}
\newcounter{tmp}

\newcommand{\A}{{\mathcal A}}

\newcommand{\C}{{\mathbb C}}
\newcommand{\D}{{\mathbb D}}

\newcommand{\K}{{\mathcal K}}

\newcommand{\R}{{\mathbb R}}

\newcommand{\es}{{\mathcal S}}

\newcommand{\CC}{{\mathcal C}}

\newcommand{\bD}{{\overline{\mathbb D}}}

\renewcommand{\Im}{{\operatorname{Im}\,}}

\renewcommand{\Re}{{\operatorname{Re}\,}}

\newcommand{\Int}{{\operatorname{Int}\,}}

\renewcommand{\arg}{\,{\operatorname{arg}\,}}

\newcommand{\Log}{{\,\operatorname{Log}\,}}
\newcommand{\Arg}{{\operatorname{Arg}\,}}
\newcommand{\aand}{{\quad\text{and}\quad}}

\begin{document}
\bibliographystyle{amsplain}
\title{
Variability regions of close-to-convex functions
}

\author[T.~Kato]{Takao Kato}
\address{Graduate School of Science and Engineering,
Yamaguchi University, Yoshida, Yamaguchi 753-8512, Japan}
\email{kato@yamaguchi-u.ac.jp}
\author[T.~Sugawa]{Toshiyuki Sugawa}
\address{Graduate School of Information Sciences,
Tohoku University, Aoba-ku, Sendai 980-8579, Japan}
\email{sugawa@math.is.tohoku.ac.jp}
\author[L.-M.~Wang]{Li-Mei Wang}
\address{School of International Technology and Management,
University of International Business and Economics, No.~10, Huixin
Dongjie, Chaoyang District, Beijing 100009, China}
\email{wangmabel@163.com}  
\begin{abstract}

M.~Biernacki gave concrete forms of the variability regions of
$z/f(z)$ and $zf'(z)/f(z)$ of close-to-convex functions $f$
for a fixed $z$ with $|z|<1$ in 1936.
The forms are, however, not necessarily convenient to determine the shape
of the full variability region of $zf'(z)/f(z)$ over all
close-to-convex functions $f$ and all points $z$ with $|z|<1.$
We will propose a couple of other forms of the variability regions
and see that the full variability region of $zf'(z)/f(z)$
is indeed the complex plane minus the origin.
We also apply them to study the variability regions of
$\log[z/f(z)]$ and $\log[zf'(z)/f(z)].$
\end{abstract}
\thanks{
The present research was supported in part by JSPS Grant-in-Aid for
Scientific Research (B) 22340025 and (C) 23540209.
}
\maketitle

\section{Introduction}
Let $\A$ denote the class of analytic functions $f$ on the unit disk
$\D=\{z\in\C: |z|<1\}$ and let $\A_0$ and $\A_1$ be its subclasses
described by the conditions $f(0)=1$ and $f(0)=f'(0)-1=0,$ respectively.
Traditionally, the subclass of $\A_1$ consisting of univalent functions
is denoted by $\es.$
A function $f$ in $\A_1$ is called {\it starlike} (resp.~{\it convex})
if $f$ is univalent and if
$f(\D)$ is starlike with respect to $0$ (resp.~convex).
It is well known that $f\in\A_1$ is starlike (resp.~convex) precisely when
$\Re [zf'(z)/f(z)]>0$ (resp.~$\Re[1+zf''(z)/f'(z)]>0$) in $|z|<1.$
The classes of starlike and convex functions in $\A_1$ will be denoted
by $\es^*$ and $\K$ respectively.

A function $f\in\A_1$ is called {\it close-to-convex} if
$\Re[e^{i\lambda}f'(z)/g'(z)]>0$ in $|z|<1$ for a convex function
$g\in\K$ and a real constant $\lambda$ with $|\lambda|<\pi/2.$
The set of close-to-convex functions in $\A_1$ will be denoted by $\CC.$
This class was first introduced and shown to be contained in $\es$
by Kaplan \cite{Kaplan52}.
A domain is called {\it close-to-convex} if it is expressed as the
image of $\D$ under the mapping $af+b$ for some $f\in\CC$ and
constants $a, b\in\C$ with $a\ne0.$
He also gave a geometric characterization in terms of turning of
the boundary of the domain.
We recommend books \cite{Duren:univ} and \cite{Goodman:univ} for general
reference on these topics.

Prior to the work of Kaplan, Biernacki \cite{Bier36} introduced the
notion of linearly accessible domains (in the strong sense).
Here, a domain in $\C$ is called linearly accessible if
its complement is a union of half-lines which do not cross each other.
Lewandowski \cite{lew58}, \cite{lew60} proved that the class of close-to-convex
domains is identical with that of linearly accessible domains
(see also \cite{BL62} and \cite{Koepf89} for simpler proofs of this fact).
Therefore, the work of Biernacki on linearly accessible domains and
their mapping functions can now be interpreted as that on close-to-convex
domains and functions.

For a non-vanishing function $g$ in $\A_0,$ unless otherwise stated,
$\log g$ will mean
the continuous branch of $\log g$ in $\D$ determined by $\log g(0)=0.$
For instance, $f(z)/z$ can be regarded as a non-vanishing function
in $\A_0$ for $f\in\es.$
Therefore, we can define $\log f(z)/z$ in the above sense.
In the present note, we are interested in the following variability regions
for a fixed $z\in\D:$
$$
\begin{array}{cclcccl}
U_z&=&\{\tfrac z{f(z)}: f\in\CC\},&
\quad&
LU_z&=&\{\log\tfrac z{f(z)}: f\in\CC\},
\\
V_z&=&\{f'(z): f\in\CC\},&
\quad&
LV_z&=&\{\log f'(z): f\in\CC\},
\\
W_z&=&\{\tfrac {zf'(z)}{f(z)}: f\in\CC\},&
\quad&
LW_z&=&\{\log\tfrac {zf'(z)}{f(z)}: f\in\CC\}.
\end{array}
$$
We collect basic properties of these sets.

\begin{lem}\label{lem:X}
\hspace{3mm}
\begin{enumerate}
\item
$X_z$ is a compact subset of $\C$ for each $z\in\D$ and
$X=U, V, W, LU, LV, LW.$
\item
$X_z=\exp(LX_z)$ for each $z\in\D$ and $X=U, V, W.$
\item
$X_z=X_r$ for $|z|=r<1$ and $X=U, V, W, LU, LV, LW.$
%For $z\in\D$ with $|z|=r,$
%$U_z=U_r, LU_z=LU_r, V_z=V_r, LV_z=LV_r, W_z=W_r, LW_z=LW_r.$
\item
$X_r\subset X_s$ for $0\le r<s<1$ and
$X=U, V, W, LU, LV, LW.$
\end{enumerate}
\end{lem}

\begin{pf}
It is enough to outline the proof since the reader can reproduce the
proof easily.
Assertion (1) follows from compactness of the family $\CC,$
whereas (2) is immediate by definition.
To see (3) and (4), it is enough to show that
$X_z\subset X_w$ for $|z|\le |w|<1.$
This follows from the fact that
the function $f_a(z)=f(az)/a$ belongs to $\CC$ again
for $f\in\CC$ and $a\in\C$ with $0<|a|\le 1.$
%The assertion $\partial X_r\cap \partial X_s=\emptyset$ for $r<s<1$
%is a consequence of the following observation due to Biernacki \cite{Bier36}:
%The extremal functions in $\CC$ for the relevant functionals consist
%of Koebe transforms of the Koebe function $z/(1-z)^2.$
\end{pf}

We remark that we can indeed show the stronger inclusion relation
$X_r\subset\Int X_s$ for $0\le r<s<1$ by observing extremal functions
corresponding to boundary points of $X_r.$
Here, $\Int E$ means the set of interior points of a subset $E$ of $\C.$
However, we do not use this property in what follows.

Set $X_{1^-}=\bigcup_{0\le r<1} X_r$ for $X=U, V, W, LU, LV, LW.$
In the sequel, $\D(a,r)$ will stand for the open disk $|z-a|<r$ in $\C$
and $\bD(a,r)$ will stand for its closure, namely, the closed disk
$|z-a|\le r.$

Biernacki \cite{Bier36} described $U_z$ and $W_z$ in his study
on linearly accessible domains and their mapping functions.
The results can be summarized as in the following.

\begin{lem}[Biernacki (1936)]\label{lem:bier}
For $0<r<1,$ the following hold:
\begin{enumerate}
\item
$U_r=\{(1+s)^2/(1+\tfrac{s+t}2): |s|\le r, |t|\le r\}
=\{2u^2/(u+v): |u-1|\le r, |v-1|\le r\}.$
\item
$W_r=(1-r^2)^{-2}U_r.$
\item\label{item:U1}
$U_{1^-}=\D(1,3)\setminus\{0\}$ and
$LU_{1^-}\subset\{w\in\C: |\Im w|<3\pi/2\}.$
\item
$LW_{1^-}\subset\{w\in\C: |\Im w|<3\pi/2\}.$
\end{enumerate}
\end{lem}

The above expressions of $U_r$ and $W_r$ are simple but somewhat implicit.
For instance, the parametrization of the boundary curve cannot be obtained immediately
and the shape of the limit $W_{1^-}$ is not clear (as we will see below,
this set is equal to $\C\setminus\{0\}$).
Therefore, it would be nice to have more explicit or more convenient expressions
of $U_r$ and $W_r$.
We propose two such expressions in the present note.

\begin{thm}\label{thm:F}
For $0<r<1,$ $U_r=F(\bD(0,r)),$ where
$$
F(z)=\frac{(3+\bar z)(1+z)^3}{3+3z+\bar z+z^2},
\quad z\in\D.
$$
\end{thm}

We will prove the theorem by describing explicitly the envelope of the family
of circles $M_s(\partial\D(0,r))$ for $s=re^{i\theta},~0\le\theta<2\pi,$
where $M_s$ is the M\"obius transformation $t\mapsto (1+s)^2/(1+(s+t)/2).$
Lewandowski \cite[p.~45]{lew60} used the envelope to prove that
the inclusion $U_r\subset\{w\in\C: \Re w\ge0\}$ (equivalently,
$W_r\subset\{w\in\C: \Re w\ge0\}$) is valid
precisely when $r\le 4\sqrt2-5.$
(This implies that the radius of starlikeness of close-to-convex functions
is $4\sqrt2-5.$)
However, any explicit form of the envelope was not given in \cite{lew60}
because it was not necessary for his results.

We note that $F(e^{i\theta})=1+3e^{i\theta}$ for $\theta\in\R,$
which agrees with Lemma \ref{lem:bier} \eqref{item:U1}.
But, this does not give enough information to determine the boundary
curve of the domain $LU_{1^-}.$
%Though an explicit form of $U_{1^-}$ is given in Lemma \ref{lem:bier},
%it is not easy to determine the shape of $LU_{1^-}.$
It turns out that $LU_{1^-}$ has relatively a simple description
though $LU_r$ does not have.
We indeed derive the following result by making use of Theorem \ref{thm:F}.

\begin{thm}\label{thm:LU1}
The variability region $LU_{1^-}$ is an unbounded Jordan domain
with the boundary curve $\gamma(t),~-2\pi<t<2\pi,$ given by
$$
\gamma(t)=
\begin{cases}
\Log(1+3e^{it}) & \quad\text{if}~ |t|< \pi \\
\Log(1-e^{it})+\dfrac{t}{|t|}\pi i & \quad\text{if}~ \pi\le|t|<2\pi.
\end{cases}
$$
\end{thm}

Here and hereafter,
$\Log w=\log|w|+i\Arg w$ denotes the principal branch of $\log w$ with
$-\pi<\Im\Log w=\Arg w\le\pi.$

As we will see in the next section, the function $\log F$ is univalent in
$\D.$ Therefore, the last theorem tells us that $F:\D\to U_{1^-}$ covers
the disk $\D(-1,1)$ bivalently whereas it covers
$\D(1,3)\setminus(\D(-1,1)\cup\{0\})$ univalently.

The following expression of $W_r$ is not very explicit but useful
in some situation.

\begin{thm}\label{thm:W}
For $0<r<1,$
$$
W_r=\left\{\frac{2u}{v(u+v)}: |u-1|\le r, |v-1|\le r\right\}.
$$
\end{thm}

Indeed,
as an application of the last theorem, we can show the following result.

\begin{thm}\label{thm:LW}
$LW_{1^-}=\{w\in\C: |\Im w|<3\pi/2\}.$
\end{thm}

Since $W_{1^-}=\exp(LW_{1^-}),$ we obtain the following corollary,
which was used in \cite{KS11PT}.

\begin{cor}
The full variability region $\{zf'(z)/f(z): z\in\D, f\in\CC\}$ is equal to
$\C\setminus\{0\}.$
\end{cor}

The corollary means that $W_{1^-}=\C\setminus\{0\}.$
We note here that this does not seem to follow immediately
from Lemma \ref{lem:bier}.

Krzy\.z \cite{Krzyz64} showed that $LV_r$ is convex and determined its shape
for $0<r<1.$

\begin{prop}[Krzy\.z]
For $0<r<1,$ the variability region $LV_r$ is convex and its boundary is
described by the curve 
$\sigma_r(t)=\log(1-re^{i\theta_2(t)})/(1-re^{i\theta_1(t)})^3,~
-\pi\le t\le \pi.$
Here,
$$
\theta_1(t)=t-\arcsin(r\sin t),\quad
\theta_2(t)=\pi+t+\arcsin(r\sin t).
$$
\end{prop}

He also proved that $LV_r$ is contained in the domain $|\Im w|<4\arcsin r$ 
for each $0<r<1$ and that this bound is sharp.
(See also \cite[Chap.~11]{Goodman:univ}.)
In particular, $LV_{1^-}\subset\{w: |\Im w|<2\pi\}.$
Since $\Re\sigma_r(t)\to+\infty$ as $r\to1^-$ for $|t|<\pi/2$
and $\Re\sigma_r(t)\to-\infty$ as $r\to1^-$ for $\pi/2<|t|\le\pi,$
it is not very easy to determine the limiting shape of $LV_r$ as $r\to1^-$
by the above proposition.
We thus complement his results by showing the following.

\begin{thm}\label{thm:v}
$$
V_r=\{(1+s)/(1+t)^3: |s|\le r, |t|\le r\}
=\{u/v^3: |u-1|\le r, |v-1|\le r\}
$$
for $0<r<1.$
Moreover, $LV_{1^-}=\{w: |\Im w|<2\pi\}$ and $V_{1^-}=\C\setminus\{0\}.$
\end{thm}

One might expect that $LU_r$ and $LW_r$ would also be convex for each
$0<r<1.$
This is, however, not true unlike $LV_r.$

\begin{thm}\label{thm:nonconvex}
The variability regions $LU_r$ and $LW_r$ are closed Jordan domains
for each $0<r<1.$
Moreover, there exists a number $0<r_0<1$ such that both $LU_r$ and $LW_r$
are not convex for every $r$ with $r_0<r<1.$
\end{thm}

We prove the above results in Section 3.
Section 2 will be devoted to the study of mapping properties of the function
$G=\log F$ that are necessary to show our results.

%\medskip
%\noindent
%\textbf{Acknowledgement.}
%The authors are grateful to the referee for corrections
%and constructive suggestions, which much improved the exposition.

\section{Univalence of the function $G=\log F$}

In order to analyze the shape of $LU_r$ or $LW_r,$ we need to
investigate mapping properties of the functions $F$ and $G=\log F,$
where $F$ is given in Theorem \ref{thm:F}.
Therefore, before showing the main results in Introduction, we see
basic properties of the functions $F$ and $G.$
Here, we remark that $F$ can be expressed in the form
$$
F(z)=\frac{(1+z)^3}{1+z\frac{3+z}{3+\bar z}}.
$$
Therefore, the continuous branch $G$ of $\log F$ with $G(0)=0$ is
represented by
\begin{equation}\label{eq:G}
G(z)=3\Log(1+z)-\Log(1+ze^{2i\phi}),\quad \phi=\Arg(3+z).
\end{equation}
The goal in this section is to prove the following:

\begin{thm}\label{thm:G}
The function $G=\log F$ is a homeomorphism of the unit disk $\D$
onto the domain $LU_{1^-}.$
\end{thm}

For $r\in(0,1)$ and $x\in(0,\pi),$ we set
$$
\Phi_r(x)=\Arg(1+re^{ix}).
$$
%More explicitly,
%$$
%\Phi_r(x)=\arcsin\left(
%\frac{r\sin x}{\sqrt{1+2r\cos x+r^2}}
%\right).
%$$
We will use the following elementary properties of the function $\Phi_r.$

\begin{lem}\label{lem:Phi}
%\begin{enumerate}
Let $r\in(0,1).$
Then
$$
\Phi_r'(x)=\frac{r(r+\cos x)}{1+2r\cos x+r^2},
\quad x\in(0,\pi).
$$
In particular, $\Phi_r(x)$ is increasing in $0<x<x_r$
and decreasing in $x_r<x<\pi,$ where
$x_r=\pi-\arccos r.$
Furthermore, $\Phi_r'(x)$ is decreasing in $0<x<\pi$
and therefore $\Phi_r$ is concave in $(0,\pi).$
%\end{enumerate}
\end{lem}

We also need the following information.

\begin{lem}\label{lem:phi}
Let $0<r<1.$
Then the inequalities $0<\theta+2\phi<\pi$ hold
for $0<\theta<\pi$ and $\phi=\Arg(3+re^{i\theta}).$
\end{lem}

\begin{pf}
Set $h_r(\theta)=\theta+2\phi
=\theta+2\Arg(3+re^{i\theta}),~0\le\theta\le\pi.$
Then
\begin{align*}
h_r'(\theta)&=1+2\frac{\partial\phi}{\partial\theta}
%=1+\frac{2r(r+3\cos\theta)}{9+6r\cos\theta+r^2}
=\frac{3(3+4r\cos\theta+r^2)}{9+6r\cos\theta+r^2} \\
&\ge\frac{3(3-4r+r^2)}{9+6r\cos\theta+r^2}
=\frac{3(1-r)(3-r)}{9+6r\cos\theta+r^2}>0.
\end{align*}
Therefore, $h_r(\theta)$ is increasing in $0<\theta<\pi,$
which implies that $0=h_r(0)<h_r(\theta)<h_r(\pi)=\pi$ for $0<\theta<\pi.$
\end{pf}

As for the function $G,$ we first show its local univalence.

\begin{lem}\label{lem:lu}
The function $G$ is orientation-preserving and locally univalent in $\D.$
\end{lem}

\begin{pf}
The partial derivatives of $G$ are given by
\begin{align*}
G_z(z)&=\frac{6+4z+3\bar z+z^2}{(1+z)(3+3z+\bar z+z^2)}, \\
G_{\bar z}(z)&=\frac{z(3+z)}{(3+\bar z)(3+3z+\bar z+z^2)}.
\end{align*}
It suffices to show that the Jacobian $J_G=|G_z|^2-|G_{\bar z}|^2$
is positive in $\D,$ which is equivalent to the condition
$|6+4z+3\bar z+z^2|>|z(1+z)|$ in $|z|<1.$
If we write $z=re^{i\theta},$ then
\begin{align*}
&\quad\ \  |6+4z+3\bar z+z^2|^2-|z(1+z)|^2 \\
&=6(6+14r\cos\theta+4r^2+6r^2\cos2\theta+r^3\cos\theta+r^3\cos3\theta) \\
&=12(1+r\cos\theta)\big(2(1+r\cos\theta)^2+1-r^2\big)>0.
\end{align*}
Thus we are done.
\end{pf}

\begin{lem}\label{lem:Re}
For a fixed $0<r<1,$ the real part of $G(re^{i\theta})$
is a decreasing function in $0\le\theta\le\pi.$
\end{lem}

\begin{pf}
By differentiating both sides of
$$
\Re G(re^{i\theta})
%=\log|F(re^{i\theta})|^2
=\frac32\log(1+2r\cos\theta+r^2)
-\frac12\log(1+2r\cos(\theta+2\phi)+r^2)
$$
with respect to $\theta,$ we obtain
$$
\frac{\partial}{\partial\theta}\Re G(re^{i\theta})
=\frac{-3r\sin\theta}{1+2r\cos\theta+r^2}
+\frac{r\sin(\theta+2\phi)}{1+2r\cos(\theta+2\phi)+r^2}
\left(1+2\frac{\partial\phi}{\partial\theta}\right).
$$
Since $\tan\phi=r\sin\theta/(3+r\cos\theta),$ we have the relations
\begin{align*}
\frac{\partial\phi}{\partial\theta}
&=\frac{r(3\cos\theta+r)}{9+6r\cos\theta+r^2}, \\
\cos 2\phi&=\frac{1-\tan^2\phi}{1+\tan^2\phi}
=1-\frac{2r^2\sin^2\theta}{9+6r\cos\theta+r^2}, \\
\sin 2\phi&=\frac{2\tan\phi}{1+\tan^2\phi}
=\frac{2r(3+r\cos\theta)\sin\theta}{9+6r\cos\theta+r^2}.
\end{align*}
We substitute them into the above expression of
$\partial(\Re G)/\partial\theta$
and make some simplifications to obtain
$$%\begin{align*}
\frac{\partial}{\partial\theta}\Re G(re^{i\theta})
=~-\frac{6r\sin\theta(3+4r\cos\theta+r^2\cos2\theta)%
H(r,\theta)}
%(9+12r\cos\theta-4r^2\sin^2\theta-4r^3\cos\theta-r^4)}%
{(1+2r\cos\theta+r^2)(9+6r\cos\theta+r^2)K(r,\theta)},
$$%\end{align*}
where
\begin{align*}
H(r,\theta)&=9+12r\cos\theta-4r^2\sin^2\theta-4r^3\cos\theta-r^4, \\
K(r,\theta)&=9-2r^2+r^4+24r\cos\theta+24r^2\cos^2\theta+8r^3\cos^3\theta.
\end{align*}
Note first that
$$
3+4r\cos\theta+r^2\cos2\theta=1-r^2+2(1+r\cos\theta)^2>0.
$$
Secondly, we compute
$$
\frac{\partial H(r,\theta)}{\partial\theta}
=-4r(3+2r\cos\theta-r^2)\sin\theta.
$$
In particular, $H(r,\theta)$ is decreasing in $0<\theta<\pi$
for a fixed $0<r<1.$
Hence,
$$
H(r,\theta)\ge H(r,\pi)=(1-r)(3-r)(3-r^2)>0.
$$
Thirdly, we see that the function $p(x)=24x+24x^2+8x^3$ is increasing
because $p'(x)=24(1+x)^2\ge0.$
Therefore,
$$
K(r,\theta)\ge K(r,\pi)=9-24r+22r^2-8r^3+r^4=(3-r)^2(1-r)^2>0.
$$
We summarize the above observations to conclude
that $\partial(\Re G)/\partial\theta<0$ for $0<\theta<\pi.$
\end{pf}

We next prove the following:

\begin{lem}\label{lem:upper}
$\Im G(z)>0$ for $z\in\D$ with $\Im z>0.$
\end{lem}

\begin{pf}
Fix $r\in(0,1).$
Let $\theta\in(0,\pi)$ and let $\phi$ be given in \eqref{eq:G}.
Note that $0<\phi<\theta.$
We need to show that
\begin{equation}\label{eq:g}
g_r(\theta)=\Im G(re^{i\theta})=3\Phi_r(\theta)-\Phi_r(\theta+2\phi)
\end{equation}
is positive.
Note that $0<\theta<\theta+2\phi<\pi$ by Lemma \ref{lem:phi}.
Lemma \ref{lem:Phi} implies that $\Phi_r$ takes its maximum value
$\Phi_r(x_r)=\arcsin r$ in $(0,\pi).$
In particular, we have $\Phi_r(\theta+2\phi)\le \arcsin r.$
Therefore, $g_r(\theta)>0$ for $x_r^-<\theta<x_r^+,$
where $x_r^-$ and $x_r^+$ are the solutions to the equation
$3\Phi_r(x)=\arcsin r$ with $0<x_r^-<x_r<x_r^+<\pi.$

We next assume that $x_r^+\le \theta<\pi.$
Since $\Phi_r$ is decreasing in $x_r^+<x<\pi$ by Lemma \ref{lem:Phi},
we see that $g_r(\theta)>\Phi_r(\theta)-\Phi_r(\theta+2\phi)>0.$
Finally, we assume that $0<\theta\le x_r^-.$
In view of concavity of $\Phi_r$ (see Lemma \ref{lem:Phi})
together with $\Phi_r(0)=0,$
we have $\Phi_r(x_r^-)=\Phi_r(x_r)/3\le\Phi_r(x_r/3).$
Hence, $x_r^-\le x_r/3.$
In particular, we have $\theta+2\phi\le 3\theta\le 3x_r^-\le x_r.$
Since $\Phi_r$ is increasing and concave in $0<x<x_r,$ the inequalities
$\Phi_r(\theta+2\phi)<\Phi_r(3\theta)\le 3\Phi_r(\theta)$ follow.
Thus we have shown that $g_r(\theta)>0$ in this case, too.
\end{pf}

The following property will be used for the proof of Theorem \ref{thm:LU1}.

\begin{lem}\label{lem:g}
The function $g_r$ defined in \eqref{eq:g} satisfies
$g_r'(\theta)>0$ for $0<\theta<x_r=\pi-\arccos r.$
\end{lem}

\begin{pf}
By definition, we have
$$
g_r'(\theta)=3\Phi_r'(\theta)-
\left(1+2\frac{\partial\phi}{\partial\theta}\right)\Phi_r'(\theta+2\phi).
$$
Since $1+2\partial\phi/\partial\theta>0$ (see the proof of Lemma \ref{lem:phi}),
Lemma \ref{lem:Phi} implies
$$
g_r'(\theta)\ge 3\Phi_r'(\theta)-
\left(1+2\frac{\partial\phi}{\partial\theta}\right)\Phi_r'(\theta)
=2\left(1-\frac{\partial\phi}{\partial\theta}\right)\Phi_r'(\theta).
$$
We note here that
$$
1-\frac{\partial\phi}{\partial\theta}
=\frac{3(3+r\cos\theta)}{9+6r\cos\theta+r^2}>0.
$$
Since $\Phi_r'(\theta)>0$ for $0<\theta<x_r$ by Lemma \ref{lem:Phi},
the required assertion follows.
\end{pf}

We are now ready to prove the theorem.

\begin{pf}[Proof of Theorem \ref{thm:G}]
Since $G$ is orientation preserving and locally univalent
by Lemma \ref{lem:lu}, it is enough to show that $G$ is injective
on the circle $|z|=r$ for each $r\in(0,1).$
We note here that $G$ is symmetric in the real axis, in other words,
$G(\bar z)=\overline{G(z)}$ for $z\in\D.$
By Lemmas \ref{lem:Re} and \ref{lem:upper}, $G$ maps the upper half of the
circle $|z|=r$ univalently onto a Jordan arc in the upper half plane.
Taking into account the symmetry, we have confirmed that
$G$ maps the circle $|z|=r$ univalently onto a Jordan curve which
is symmetric in the real axis.
Thus the proof is complete.
\end{pf}

\section{Proof of main results}

In this section, we show the main results presented in Section 1.
We begin with the proof of Theorem \ref{thm:F}.

\begin{pf}[Proof of Theorem \ref{thm:F}]

Fix $0<r<1.$
To simplify computations, we consider the set $1/U_r=\{1/w: w\in U_r\}$
instead of $U_r.$
For $s\in\D,$ we denote by $\Delta_s$ the image of
the disk $|t|\le r$ under the mapping 
$t\mapsto(1-\frac{s+t}2)/(1-s)^2.$
Note that $\Delta_s$ is the closed disk with center
$(1-s/2)/(1-s)^2$ and radius $r|1-s|^{-2}/2.$
It is easily verified that the function $(1-s/2)/(1-s)^2$ is
univalent on the unit disk $|s|<1.$
For $0\le q\le r,$
let $E_q$ be the union of $\Delta_s$ over $s\in\C$
with $|s|=q.$
Then Lemma \ref{lem:bier} (1) implies that 
$E_q$ sweeps $1/U_r$ when $q$ moves from $0$ to $r.$
Note here that $E_0$ is the disk $\bD(1,r/2).$
By these observations, we see that the boundary of $1/U_r$
is contained in the {\it outer} envelope of the family of
circles $\partial\Delta_s$ over $|s|=r.$

We next describe the outer envelope.
Let $c(\alpha)$ and $\rho(\alpha)$ be the center and the radius
of the disk $\Delta_s$ for $s=re^{i\alpha}.$
Explicitly,
$$
c(\alpha)=\frac{2-s}{2(1-s)^2}
\aand
\rho(\alpha)=\frac{r}{2|1-s|^2}=\frac{r}{2(1-2r\cos\alpha+r^2)}.
$$
Note that
$$
c'(\alpha)=\frac{is(3-s)}{2(1-s)^3}
\aand
\rho'(\alpha)=-\frac{r^2\sin\alpha}{|1-s|^4}.
$$
We can parametrize the outer envelope in the form
$$
\zeta(\alpha)=c(\alpha)+\rho(\alpha)e^{i\beta(\alpha)},\quad
-\pi<\alpha\le\pi.
$$
By the symmetry in the real axis and the fact that $\Im c'(0)>0,$
we can take $\beta(\alpha)$ so that $\beta(0)=0.$
Here, $\beta=\beta(\alpha)$ is a real-valued differentiable function 
of $\alpha$ satisfying the condition that the tangent vector
$\zeta'(\alpha)$ is tangent to the circle
$|w-c(\alpha)|=\rho(\alpha)$ at $\zeta(\alpha).$
In other words, $\zeta'(\alpha)=kie^{i\beta}$ for a real number $k.$
Taking the real part of the relation
$$
\zeta'(\alpha)e^{-i\beta}
=c'(\alpha)e^{-i\beta}+\rho'(\alpha)+i\beta'(\alpha)\rho(\alpha)=ki,
$$
we obtain
$$
\Re[c'(\alpha)e^{-i\beta}]+\rho'(\alpha)=0,
$$
which implies that
$$
\cos(\beta-\arg c'(\alpha))=
\cos(\arg c'(\alpha)-\beta)=-\frac{\rho'(\alpha)}{|c'(\alpha)|}
=\frac{2r\sin\alpha}{|1-s||3-s|}.
$$
Hence,
\begin{align*}
\frac{|c'(\alpha)|}{c'(\alpha)}e^{i\beta}=
e^{i(\beta-\arg c'(\alpha))}
&=\frac{2r\sin\alpha}{|1-s||3-s|}
\pm i\sqrt{1-\left(\frac{2r\sin\alpha}{|1-s||3-s|}\right)^2} \\
&=\frac{2r\sin\alpha\pm i(3-4r\cos\alpha+r^2)}{|1-s||3-s|}.
\end{align*}

We recall that $\beta(0)=0$ and substitute $\alpha=0$
into this relation in order to eliminate ambiguity of the sign.
We then see that the minus sign should be taken there.
Hence,
$$
e^{i\beta}=\frac{-i(1-\bar s)(3-s)}{|1-s||3-s|}
\cdot\frac{c'(\alpha)}{|c'(\alpha)|}
=\frac{s(3-s)(1-\bar s)^2}{r(3-\bar s)(1-s)^2}.
$$
We now get the form of $\zeta(\alpha):$
\begin{align*}
\zeta(\alpha)&=\frac{2-s}{2(1-s)^2}+\frac{r}{2|1-s|^2}\cdot
\frac{s(3-s)(1-\bar s)^2}{r(3-\bar s)(1-s)^2} \\
&=\frac{3-3s-\bar s+s^2}{(3-\bar s)(1-s)^3}=\frac1{F(-s)}.
\end{align*}
Thus the assertion follows.
\end{pf}

\begin{pf}[Proof of Theorem \ref{thm:LU1}]
As we saw, the mapping $G=\log F$ is a homeomorphism of $\D$
onto $LU_{1^-}.$
We now observe that $G$ extends continuously to $\overline\D\setminus\{-1\}.$
Since $G(e^{it})=\gamma(t)$ for $|t|<\pi,$ the boundary of $LU_{1^-}$
contains the arc $\gamma([-\pi,\pi]).$

We next investigate the limit points of $G(z)$ as $z\to -1.$
Let $\alpha(\delta)=(a\delta)^{1/3}$ and put
$z=(1-\delta)e^{i(\pi-\alpha(\delta))}$ for $0<\delta<1$
and a positive constant $a.$
Then
\begin{align*}
z&=-(1-\delta)\left(1-i\alpha(\delta)-\frac{\alpha(\delta)^2}2
+\frac{i\alpha(\delta)^3}6+O(\delta^{4/3})\right)\\
&=-1+i(a\delta)^{1/3}+\frac{(a\delta)^{2/3}}2
+\left(1-\frac{ia}6\right)\delta+O(\delta^{4/3})
\end{align*}
as $\delta\to0^+.$
Therefore,
$$
\frac{3+z}{3+\bar z}
=1+i(a\delta)^{1/3}-\frac{(a\delta)^{2/3}}2-\frac{2ia\delta}3+O(\delta^{4/3})
$$
and
$$
1+z\frac{3+z}{3+\bar z}=\left(1+\frac{ia}2\right)\delta+O(\delta^{4/3}).
$$
Since $(1+z)^3=-ia\delta+O(\delta^{4/3}),$ we have
$$
\frac{(1+z)^3}{1+z\frac{3+z}{3+\bar z}}
=-\left(1+\frac{a+2i}{a-2i}\right)+O(\delta^{1/3})
$$
as $\delta\to0^+.$
Thus
$$
\lim_{\delta\to0^+}G((1-\delta)e^{i(\pi-\alpha(\delta))})
=\pi i+\log\left(1+\frac{a+2i}{a-2i}\right).
$$
Since $a$ is an arbitrary positive real number, the boundary of $LU_{1^-}=G(\D)$
contains the curve $\gamma(t):~\pi<t<2\pi.$
The same is true for $-2\pi<t<-\pi$ by the symmetry of the function $G.$

The remaining thing is to prove that the boundary of $LU_{1^-}$
in $\C$ contains no other points than the curve $\Gamma=\{\gamma(t):|t|<2\pi\}.$
We note here that $LU_r$ is convex in the direction of imaginary axis
for each $0<r<1$ by Lemma \ref{lem:Re}.
Therefore, the same is true for the limit $LU_{1^-}.$
We observe also that $\Gamma$ encloses an unbounded Jordan domain
convex in the direction of imaginary axis.

Suppose that there is a boundary point $p_0$ of $LU_{1^-}$ with
$p_0\notin\Gamma.$
We may assume that $\Im p_0>0.$
Let $p_1$ be the point in $\Gamma$ with $\Im p_1>0$ and $\Re p_1=\Re p_0.$
Then the convexity of $LU_{1^-}$ in the direction of imaginary axis implies
that the segment $[p_0,p_1]$ is contained in $\partial LU_{1^-}.$
We can choose $p_0$ so that the segment is maximal.
Since the family of smooth Jordan domains $LU_r,~0<r<1,$ exhausts the domain
$LU_{1^-},$ for a small enough $\delta>0$
there exist three points $z_1^-(\delta), z_0(\delta), z_1^+(\delta)$
on the circle $|z|=1-\delta$ with
$0<\Arg z_1^-(\delta)<\Arg  z_0(\delta)<\Arg z_1^+(\delta)$ such that
$G(z_1^-(\delta))\to p_1,~G(z_0(\delta))\to p_0,~G(z_1^+(\delta))\to p_1$ as
$\delta\to0^+.$
In particular, $\Im G(z_0(\delta))<\Im G(z_1^\pm(\delta))$ for sufficiently
small $\delta>0.$
Therefore, $g_{1-\delta}(\theta)=\Im G((1-\delta)e^{i\theta})$
takes a local minimum at a point
$\theta_0$ with $\Arg z_1^-(\delta)<\theta_0<\Arg z_1^+(\delta).$
In particular, $g_{1-\delta}'(\theta_0)=0.$
Note here that
$$
\Re G((1-\delta)e^{i\theta_0})\to \Re p_0\quad(\delta\to0^+).
$$
We write $\theta_0=\pi-\beta(\delta).$
Then, by Lemma \ref{lem:g}, we see that $\theta_0\ge x_{1-\delta},$
equivalently, $\beta(\delta)\le \arccos(1-\delta).$
This implies that $\beta(\delta)=O(\delta^{1/2})$ as $\delta\to0^+.$
Therefore, $z=(1-\delta)e^{i(\pi-\beta(\delta))}=-1+i\beta(\delta)+O(\delta),$
$(3+z)/(3+\bar z)=1+i\beta(\delta)+O(\delta)$
and thus $1+z(3+z)/(3+\bar z)=O(\delta)$ as $\delta\to0^+.$
In particular,
$$
\Re G((1-\delta)e^{i\theta_0})\to -\infty\quad(\delta\to0^+),
$$
which is a contradiction.

We now conclude that $\partial LU_{1^-}=\Gamma.$
\end{pf}

In order to prove Theorem \ref{thm:W}, we will make use of a weakened
version of Lemma 5.1 in Greiner and Roth \cite{GR03}, which is an outcome
of the duality methods developed by Ruscheweyh and Sheil-Small.

For $|a|\le1, |b|\le 1,$ define a function $f_{a,b}\in\A_1$ by
$$
f_{a,b}(z)=z\frac{1+(a+b)z/2}{(1+bz)^2}.
$$
It is easy to see that $f_{a,b}$ belongs to the class $\CC$
of close-to-convex functions.
The linear space $\A$ is naturally equipped with the topology of
uniform convergence on compact subsets in $\D.$

\begin{lem}[$\text{\cite[Lemma 5.1]{GR03}}$]\label{lem:GR}
Let $\lambda_1$ and $\lambda_2$ be continuous linear functionals
on $\A$ such that $\lambda_2$ does not vanish on $\CC.$
Then for every $f\in\CC$ there exist complex numbers $a,b$
with $|a|\le 1, |b|\le 1$ such that
$$
\frac{\lambda_1(f)}{\lambda_2(f)}
=\frac{\lambda_1(f_{a,b})}{\lambda_2(f_{a,b})}.
$$
\end{lem}

\begin{pf}[Proof of Theorem \ref{thm:W}]
Fix $0<r<1.$
Let $f\in\CC.$
We now apply Lemma \ref{lem:GR} to the choice
$\lambda_1(f)=rf'(r)$ and $\lambda_2(f)=f(r)$ to see that
$$
\frac{rf'(r)}{f(r)}=\frac{\lambda_1(f)}{\lambda_2(f)}
=\frac{\lambda_1(f_{a,b})}{\lambda_2(f_{a,b})}
=\frac{2(1+ar)}{(1+br)(2+(a+b)r)}
$$
for some $a,b\in\overline\D.$
The proof is complete by letting $u=1+ar$ and $v=1+br$.
\end{pf}

\begin{pf}[Proof of Theorem \ref{thm:LW}]
Let $\Omega=\{(r,s,t)\in\R^3: 0<s<2, 0<rs^2<2, -\pi/2<t<\pi/2\}.$
Then $u=rs^2e^{it}\cos^2t$ and $v=se^{-it}\cos t$ satisfy
$|u-1|<1$ and $|v-1|<1,$ whence the point
$$
w(r,s,t)=\log\frac{2u}{v(u+v)}=\log (2r)+3it-\log (1+rse^{2it}\cos t)
$$
belongs to the region $LW_{1^-}$ for $(r,s,t)\in\Omega$
by Theorem \ref{thm:W}.

For a given point $z_0=x_0+iy_0$ with $|y_0|<3\pi/2,$
we now look for $(r,s,t)\in\Omega$ such that $w(r,s,t)=z_0.$

Let $r_0=e^{x_0}/2$ and take small enough $0<s_0<1$ so that
$r_0s_0<1/2.$
Then $r_0s_0^2<s_0<2$ and $x_0\pm3\pi i/2$
are the endpoints of the curve $\alpha(t)=w(r_0,s_0,t),~ -\pi/2<t<\pi/2.$
We now take a $t_0\in(-\pi/2,\pi/2)$ such that $\Im\alpha(t_0)=y_0$
and let $x_1=x_0-\Re\alpha(t_0).$
Since the function $-\log(1-x)$ is convex, we have the inequality
$-\log(1-x)\le 2x\log 2$ for $0\le x<1/2.$
We now estimate $-x_1$ in the following way:
$$
-x_1=-\log|1+r_0s_0e^{2it_0}\cos t_0|
\le -\log(1-r_0s_0)
\le 2r_0s_0\log 2,
$$
which implies
$$
r_0s_0^2e^{-x_1}<s_0e^{-x_1}\le s_0e^{2r_0s_0\log 2}
<s_0e^{\log 2}=2s_0<2.
$$
Therefore $(r_0e^{x_1},s_0e^{-x_1},t_0)\in\Omega$ and
$$
w(r_0e^{x_1},s_0e^{-x_1},t_0)=x_1+w(r_0,s_0,t_0)=x_0+iy_0=z_0
$$
as desired.
\end{pf}

\begin{pf}[Proof of Theorem \ref{thm:v}]
For a fixed $0<r<1,$ we consider the continuous linear functionals
$\lambda_1$ and $\lambda_2$ on $\A$ defined by
$\lambda_1(f)=f'(r)$ and $\lambda_2(f)=f'(0)$ for $f\in\A.$
Then Lemma \ref{lem:GR} implies that for any $f\in\CC,$
$$
f'(r)=\frac{\lambda_1(f)}{\lambda_2(f)}
=\frac{\lambda_1(f_{a,b})}{\lambda_2(f_{a,b})}
=\frac{1+ar}{(1+br)^3}
$$
for some $a,b\in\overline{\D}.$
Thus the first part of the theorem has been proved.

By the first part, we have the expression
$LV_{1^-}=\{\Log(1+z)-3\Log(1+w): z,w\in\D\}.$
Let $a$ and $b$ be real numbers with $|b|<\pi/2.$
We shall show that $a+4bi\in LV_{1^-}.$
It is easy to observe that the domain $\{\Log(1+z): z\in\D\}$
is convex and its boundary curve 
$$
\tau(t)=\Log(1+e^{it})=\log(2\cos\tfrac t2)+\tfrac t2 i
\qquad (-\pi<t<\pi)
$$
satisifies
$\Re\tau(t)\to -\infty$ and $\Im\tau(t)\to\pm\pi/2$ as $t\to\pm\pi^\mp.$
Therefore, there are $z,w\in\D$ such that
$a-3c+bi=\Log(1+z)$ and $-c-bi=\Log(1+w)$ for a sufficiently large $c>0.$
In particular, $a+4bi=\Log(1+z)-3\Log(1+w)\in LV_{1^-}.$
\end{pf}

\begin{pf}[Proof of Theorem \ref{thm:nonconvex}]
Since $LU_r=LW_r+\log(1-r^2)$ by Lemma \ref{lem:bier} (2), 
it suffices to prove the assertion for $LU_r.$
If there is no such an $r_0$ as in the assertion, then
the limiting domain $LU_{1^-}$ must be convex.
Note that $LU_{1^-}$ is convex if and only if
$\frac{d}{dt}\arg \gamma'(t)\ge0,$
where $\gamma$ is given in Theorem \ref{thm:LU1}.
A simple computation gives us
\begin{align*}
\frac{d}{dt}\arg \gamma'(t)&=\Im\frac{d}{dt}\log\gamma'(t) \\
&=\Re \frac1{1+3e^{it}}=\frac{1+3\cos t}{|1+3e^{it}|^2}
\end{align*}
for $|t|<\pi.$
This is negative when $\cos t<-1/3$ and thus we get a contradiction.
The proof is now complete.
\end{pf}

\def\cprime{$'$} \def\cprime{$'$} \def\cprime{$'$}
\providecommand{\bysame}{\leavevmode\hbox to3em{\hrulefill}\thinspace}
\providecommand{\MR}{\relax\ifhmode\unskip\space\fi MR }
\providecommand{\MRhref}[2]{%
  \href{http://www.ams.org/mathscinet-getitem?mr=#1}{#2}
}
\providecommand{\href}[2]{#2}

%\bibliography{papers}
\end{document}